# Properties and Problems Related to the Smarandache Type Functions

*Sebastian Martin Ruiz, Avda. De Regla, 43, Chipiona 11550 (Cadiz), Spain*
*M. Perez, Rehoboth, Box 141, NM 87301, USA*

**Abstract**:
In this paper we present the definitions and some properties of several Smarandache type functions that are involved in many proposed solved and unsolved problems and conjectures in number theory and recreational mathematics. Examples are also provided. Interesting solved problems related to them are attached as addenda to this article.

**1. Introduction**:
The most known *Smarandache Function*, which became a classical function in the number theory, is the following:

S: N* → N*, S(1) = 1, S(n) is the smallest integer such that S(n)! is divisible by n.

For example: S(6)=3, because 3! = 1·2·3 = 6, which is divisible by 6, and 3 is the smallest number with this property, i.e. 2! is not divisible by 6.
S(8) = 4, S(11) = 11.
This function has been very much studied in the last decade and interesting properties have been found related to it.

**2. Properties**:
2.1. Max {p, p prime and p divides n} ≤ S(n) ≤ n for any positive integer n.
2.2. If $n = (p_1\text{^}s_1) \cdot (p_2\text{^}s_2) \cdot\ldots\cdot (p_k\text{^}s_k)$, where $p_1, p_2, \ldots, p_k$ are distinct primes, then $S(n) = \max_{1\le i\le k} \{S(p_i\text{^}s_i)\} \le \max_{1\le i\le k} \{ p_i \cdot s_i \}$.

2.3. *Characterization of a prime number*:

Let p be an integer > 4. Then: p is prime if and only if S(p) = p.

Proof:
Let p be prime > 4, and suppose S(p) = m < p, then m! is not be divisible by p, therefore S(p) = p.
Now, let S(p) = p and p ≠ 4; suppose p is not prime, whence there exist two integers s and t, with s ≤ t < p, such that p = s·t, but then S(p) ≤ t ≠ p because t! is divisible by s and by t in the same time (i.e. t! is divisible by p). Contradiction.

2.4 *An exact formula to calculate the number of primes less than or equal to x* (L. Seagull):
If x is an integer ≥ 4, then the number of prime numbers ≤ x is:

$$\Pi(x) = -1 + \sum_{k=2}^{x} \left\lfloor \frac{S(k)}{k} \right\rfloor$$

where S(k) is the classical Smarandache Function, and ⌊a⌋
means the interior integer part of a (the smallest integer greater than or equal to a).

Proof:
Knowing the Smarandache Function has the property that if p > 4 then S(p) = p if only if p is prime, and S(k) <= k for any k, and S(4) = 4 (the only exception from the first rule), we easily find an exact formula for the number of primes less than or equal to x.

**4. Conjectures**:
4.1. The diophantine equation S(n) = S(n+1) has no solutions. (L. Tutescu)
4.2. The diophantine equation S(n) + S(n+1) = S(n+2) has infinitely many solutions. (I. M. Radu)

**5. More Smarandache Type Functions** have been also considered and studied, such as:

5.1. *Smarandache Double Factorial Function*, Sdf(n) is the smallest integer such that Sdf(n)!! is divisible by n, where the double factorial m!! = 1·3·5·...·m, if m is odd;
and m!! = 2·4·6·...·m, if m is even.

For example:

| | | | | | | | | | | | | | | | | |
|---|---|---|---|---|---|---|---|---|---|---|---|---|---|---|---|---|
| n | 1 | 2 | 3 | 4 | 5 | 6 | 7 | 8 | 9 | 10 | 11 | 12 | 13 | 14 | 15 | 16 |
| SDF(n) | 1 | 2 | 3 | 4 | 5 | 6 | 7 | 4 | 9 | 10 | 11 | 6 | 13 | 14 | 5 | 6 |

5.2. *Smarandache-Kurepa Function*:
For p prime, SK(p) is the smallest integer such that !SK(p) is divisible by p, where !m = 0! + 1! + 2! + ... + (m-1)!

For example:

| | | | | | | | | | | | | | | |
|---|---|---|---|---|---|---|---|---|---|---|---|---|---|---|
| p | 2 | 3 | 7 | 11 | 17 | 19 | 23 | 31 | 37 | 41 | 61 | 71 | 73 | 89 |
| SK(p) | 2 | 4 | 6 | 6 | 5 | 7 | 7 | 12 | 22 | 16 | 55 | 54 | 42 | 24 |

5.3. *Smarandache-Wagstaff Function*:
For p prime, SW(p) is the smallest integer such that W(SW(p)) is divisible by p, where W(m) = 1! + 2! + ... + (m)!

For example:

| | | | | | | | | | | | | | |
|---|---|---|---|---|---|---|---|---|---|---|---|---|---|
| p | 3 | 11 | 17 | 23 | 29 | 37 | 41 | 43 | 53 | 67 | 73 | 79 | 97 |
| SW(p) | 2 | 4 | 5 | 12 | 19 | 24 | 32 | 19 | 20 | 20 | 7 | 57 | 6 |

5.4. *Smarandache Ceil Functions of k-th Order*:
Sk(n) is the smallest integer for which n divides Sk(n)^k.

For example, for k=2, we have:

| | | | | | | | | | | | | | | | | |
|---|---|---|---|---|---|---|---|---|---|---|---|---|---|---|---|---|
| n | 1 | 2 | 3 | 4 | 5 | 6 | 7 | 8 | 9 | 10 | 11 | 12 | 13 | 14 | 15 | 16 |
| S2(n) | 1 | 2 | 3 | 2 | 5 | 6 | 7 | 4 | 3 | 10 | 11 | 6 | 13 | 14 | 15 | 4 |

5.5. *Pseudo-Smarandache Function*:

Z(n) is the smallest integer such that 1 + 2 + ... + Z(n) is divisible by n.

For example:

| n | 1 | 2 | 3 | 4 | 5 | 6 | 7 |
|---|---|---|---|---|---|---|---|
| Z(n) | 1 | 3 | 2 | 3 | 4 | 3 | 6 |

5.6. *Smarandache Near-To-Primordial Function*:

SNTP(n) is the smallest prime such that either p# - 1, p# , or p# + 1 is divisible by n, where p#, of a prime number p, is the product of all primes less than or equal to p.

For example:

| n | 1 | 2 | 3 | 4 | 5 | 6 | 7 | 8 | 9 | 10 | 11 | ... | 59 | ... |
|---|---|---|---|---|---|---|---|---|---|---|---|---|---|---|
| SNTP(n) | 2 | 2 | 2 | 5 | 3 | 3 | 3 | 5 | ? | 5 | 11 | ... | 13 | ... |

## 6. Other Smarandache Type Functions also studied in the last years are:

6.1. Let f: Z ↦ Z be a strictly increasing function and x an element in R. Then:

a) *Inferior Smarandache f-Part* of x, ISf(x), is the smallest k such that f(k) ≤ x < f(k+1).

b) *Superior Smarandache f-Part* of x, SSf(x) is the smallest k such that f(k) < x ≤ f(k+1).

Particular cases:

a) Inferior S-Prime Part:

For any positive real number n one defines ISp(n) as the largest prime number less than or equal to n.

The first values of this function are:

2,3,3,5,5,7,7,7,7,11,11,13,13,13,13,17,17,19,19,19,19,23,23.

b) Superior S-Prime Part:

For any positive real number n one defines SSp(n) as the smallest prime number greater than or equal to n.

The first values of this function are:

2,2,2,3,5,5,7,7,11,11,11,11,13,13,17,17,17,17,19,19,23,23,23.

c) Inferior S-Square Part:

For any positive real number n one defines ISs(n) as the largest square less than or equal to n.

The first values of this function are:

0,1,1,1,4,4,4,4,4,9,9,9,9,9,9,9,16,16,16,16,16,16,16,16,16,25,25.

b) Superior S-Square Part:

For any positive real number n one defines SSs(n) as the smallest square greater than or equal to n.

The first values of this function are:

0,1,4,4,4,9,9,9,9,9,16,16,16,16,16,16,16,25,25,25,25,25,25,25,25,25,36.

d) Inferior S-Cubic Part:

For any positive real number n one defines ISc(n) as the largest cube less than or equal to n.

The first values of this function are:
0,1,1,1,1,1,1,1,8,8,8,8,8,8,8,8,8,8,8,8,8,8,8,8,8,8,27,27,27,27.

e) Superior S-Cube Part:
For any positive real number n one defines SSs(n) as the smallest cube greater than or equal to n.
The first values of this function are:
0,1,8,8,8,8,8,8,8,27,27,27,27,27,27,27,27,27,27,27,27,27,27,27,27,27.

f) Inferior S-Factorial Part:
For any positive real number n one defines ISf(n) as the largest factorial less than or equal to n.
The first values of this function are:
1,2,2,2,2,6,6,6,6,6,6,6,6,6,6,6,6,6,6,6,6,6,24,24,24,24,24,24,24.

g) Superior S-Factorial Part:
For any positive real number n one defines SSf(n) as the smallest factorial greater than or equal to n.
The first values of this function are:
1,2,6,6,6,6,24,24,24,24,24,24,24,24,24,24,24,24,24,24,24,24,24,24,120.

Remark 1: This is a generalization of the inferior/superior integer part of a number.

6.2. Let f: Z ↦ Z be a strictly increasing function and x an element in R. Then:
*Fractional Smarandache f-Part* of x, FSf(x) = x - ISf(x), where ISf(x) is the Inferior Smarandache f-Part of x defined above.

Particular cases:

a) Fractional S-Prime Part:
FSp(x) = x - ISp(x),
where ISp(x) is the Inferior S-Prime Part defined above.
Example: FSp(12.501) = 12.501 - 11 = 1.501.

b) Fractional S-Square Part:
FSs(x) = x - ISs(x),
where ISs(x) is the Inferior S-Square Part defined above.
Example: FSs(12.501) = 12.501 - 9 = 3.501.

c) Fractional S-Cubic Part:
FSc(x) = x - ISc(x),
where ISc(x) is the Inferior S-Cubic Part defined above.
Example: FSc(12.501) = 12.501 - 8 = 4.501.

d) Fractional S-Factorial Part:
FSf(x) = x - ISf(x),
where ISf(x) is the Inferior S-Factorial Part defined above.
Example: FSf(12.501) = 12.501 - 6 = 6.501.

Remark 2.1: This is a generalization of the fractional part of a number.

Remark 2.2: In a similar way one defines:
- the Inferior Fractional Smarandache f-Part:
IFSf(x) = x - ISf(x) = FSf(x);
- and the Superior Fractional Smarandache f-Part:
SFSf(x) = SSf(x) - x;
for example: Superior Fractional S-Cubic Part of 12.501

$= 27 - 12.501 = 14.499.$

6.3. Let g: A $\mapsto$ A be a strictly increasing function, and let "~" be a given internal law on A. Then we say that f: A $\mapsto$ A is *smarandachely complementary with respect to the function g and the internal law "~"* if:

f(x) is the smallest k such that there exists a z in A so that $x \sim k = g(z)$.

Particular cases:

a) S-Square Complementary Function:
f: N $\mapsto$ N, f(x) = the smallest k such that $x \cdot k$ is a perfect square.
The first values of this function are:
1,2,3,1,5,6,7,2,1,10,11,3,14,15,1,17,2,19,5,21,22,23,6,1,26,3,7.

b) S-Cubic Complementary Function:
f: N $\mapsto$ N, f(x) = the smallest k such that $x \cdot k$ is a perfect cube.
The first values of this function are:
1,4,9,2,25,36,49,1,3,100,121,18,169,196,225,4,289,12,361,50.

More generally:

c) S-m-power Complementary Function:
f: N $\mapsto$ N, f(x) = the smallest k such that $x \cdot k$ is a perfect m-power.

d) S-Prime Complementary Function:
f: N $\mapsto$ N, f(x) = the smallest k such that x+k is a prime.
The first values of this function are:
1,0,0,1,0,1,0,3,2,1,0,1,0,3,2,1,0,1,0,3,2,1,0,5,4,3,2,1,0,1,0,5.

6.4. *S-Multiplicative Function*:
is a function $f : N^* \mapsto N^*$ such that for any (a, b) = 1, $f(a \cdot b) = \max \{f(a), f(b)\}$;
[i.e. it reflects the main property of the Smarandache function].
The following functions are obviously S-multiplicative:
a) The constant function $f : N^* \mapsto N^*$, f(n) = 1.
b) The Smarandache function S: $N^* \mapsto N$, $S(n) = \max\{ p \mid p! : n\}$.

Certainly, many properties of multiplicative functions can be translated for S-multiplicative functions.

## 7. Functional Smarandache Iterations:

7.1. *Functional Smarandache Iteration of First Kind*:

Let f: A $\mapsto$ A be a function, such that $f(x) \leq x$ for all x, and
$\min \{f(x), x \in A\} \geq m_0 \neq -\infty$.

Let f have $p \geq 1$ fix points: $m_0 \leq x_1 < x_2 < \ldots < x_p$.
[The point x is called fix if f(x) = x.] Then:

$SI1_f(x)$ = the smallest number of iterations k such that
f(f(...f(x)...)) = constant.

iterated k times

Example:

Let $n > 1$ be an integer, and $d(n)$ be the number of positive divisors of n,
$d: N \mapsto N$.
Then $SI1_d(n)$ is the smallest number of iterations k
such that $d(d(...d(n)...)) = 2$;
iterated k times

because $d(n) < n$ for $n > 2$, and the fix points of the function d are 1 and 2.
Thus $SI1_d(6) = 3$, because $d(d(d(6))) = d(d(4)) = d(3) = 2$ = constant.
$SI1_d(5) = 1$, because $d(5) = 2$.

7.2. *Functional Smarandache Iteration of Second Kind*:

Let $g: A \mapsto A$ be a function, such that $g(x) > x$ for all x, and let $b > x$. Then:
$SI2_g(x, b)$ = the smallest number of iterations k such that
$g(g(...g(x)...)) \geq b$.
iterated k times

Example:

Let $n > 1$ be an integer, and $\Sigma(n)$ be the sum of positive divisors of n (1 and n included), $\Sigma: N \mapsto N$.

Then $SI2_\Sigma(n, b)$ is the smallest number of iterations k such that
$\Sigma(\Sigma(...\Sigma(n)...)) \geq b$,
iterated k times
because $\Sigma(n) > n$ for $n > 1$.
Thus $SI2_\Sigma(4, 11) = 3$, because $\Sigma(\Sigma(\Sigma(4))) = \Sigma(\Sigma(7)) = \Sigma(8) = 15 \geq 11$.

7.3. *Functional Smarandache Iteration of Third Kind*:

Let $h: \mapsto A$ be a function, such that $h(x) < x$ for all x, and let $b < x$. Then:
$SI3_h(x, b)$ = the smallest number of iterations k such that
$h(h(...h(x)...)) \leq b$.
iterated k times

Example:

Let n be an integer and gd(n) be the greatest divisor of n, less than n,
$gd: N^* \mapsto N^*$. Then $gd(n) < n$ for $n > 1$.

$SI3_{gd}(60, 3) = 4$, because $gd(gd(gd(gd(60)))) = gd(gd(gd(30))) = gd(gd(15)) = gd(5) = 1 \leq 3$.

Other e-books on Smarandache Functions and Sequences can be downloaded from http://www.gallup.unm.edu/~smarandache/eBooks-otherformats.htm
and research papers from http://www.gallup.unm.edu/~smarandache/math.htm.